\documentclass[reqno]{amsproc}
\usepackage{amssymb}
\usepackage{amsmath}
\usepackage{mathdots}
\usepackage{amsbsy}
\usepackage{amscd}
\usepackage{amsthm}
\renewcommand{\thefootnote}{}

\usepackage{url}
\usepackage{amsmath,amssymb,amsthm,amsfonts,longtable}
\usepackage[colorlinks,citecolor=red,urlcolor=blue,bookmarks=false,hypertexnames=true]{hyperref} 

\newcommand{\Mod}[1]{\ (\mathrm{mod}\ #1)}

\begin{document}

\title{Presentations for the groups of order $p^6$ for prime $p\geq 7$}
\author{M.F.\ Newman, E.A. O'Brien and M.R.\ Vaughan-Lee}

\begin{abstract}
We list those presentations for the groups of order 
$p^6$ for prime $p\geq 7$ which are 
distributed as part of the {\sc SmallGroups} library. 
\end{abstract}

\maketitle 
\footnote{
We thank Sunil Kumar Prajapati and Ayush Udeep
for their assistance in preparing and checking the content
of Table~\ref{present}.  This work was supported in part by the Marsden Fund 
of New Zealand  via grant 20-UOA-107.

}

In \cite{p6} we proved that there are
\begin{center}
$3p^2+39p+344+24 \gcd(p-1,3)+11 \gcd(p-1,4)+2\gcd(p-1,5)$
\end{center}
isomorphism types of groups of order $p^6$
for every prime $p \ge 7$.

Our determination of these groups differs significantly from that claimed
by James \cite{James80}; some of the errors
in his work are recorded in \cite{error-list}.
See \cite{p6} for details of checks performed on our determination. 
We refer the reader interested in the chequered history of 
classification of these groups to \cite{p6}.

Our presentations for these groups are available 
publicly as part of the {\sc SmallGroups} library \cite{Small}. 
There continues to be demand for a human-readable description of the groups.
One possible approach is to apply the Baker-Campbell-Hausdorff formula to 
the corresponding Lie ring presentations listed in \cite{mrvl}.
Here we provide group presentations explicitly: 
the parametrised presentations listed in Table \ref{present} mirror 
those stored in the realisations of {\sc SmallGroups} distributed in 
{\sf GAP} \cite{GAP4} and {\sc Magma} \cite{Magma}.
An electronic version of the latter code 
is available at \cite{github}.

The groups are classified into 43 isoclinism families labelled 
$\Phi_{i}$ for $1 \leq i \leq 43$.
Invariants for each family are tabulated in \cite{James80};
these include central quotient, derived group, lower central 
series structure, conjugacy class structure, and 
number of irreducible characters.

The ordering of the rows in Table \ref{present} is based on the ordering of the 
isoclinism families
to which the groups belong. We follow the ordering of these families given by 
Easterfield \cite{Easterfield40}. This is first by rank. 
The ordering within families follows that 
of \cite{Easterfield40} -- except for $\Phi_{19}$ which follows 
\cite{James80}. 

The families of rank 6 are $\Phi_{11}$ to $\Phi_{43}$.
These are ordered by decreasing size of the centres. 
With the exception of $\Phi_{19}$, 
each listed presentation has generating set 
$\{\alpha_1,\ldots,\alpha_6\}$ taken from a composition series 
of the group which refines its upper central series,
so $\alpha_1$ is always central. 
The relation set consists of power relations and fifteen commutator 
relations $[\alpha_i,\alpha_j] = \ldots$ 
-- we omit those commutator relations which are trivial.
Some of the relations involve $\nu$, the smallest positive integer which 
is a quadratic non-residue mod $p$,
and some $\omega$, the smallest positive integer which is 
a primitive root mod $p$.
For some rows there are one or two additional parameters.
Examples are: 
\begin{itemize}
\item 
$G_{(11,14r)}$ where $r$ runs over ${1,\nu}$; 
\item 
$G_{(16,12r)}$ where $r=1,\omega,\omega^2$ when $p\equiv 1\bmod 3$, and $r=1$ when $p\equiv 2\bmod 3$; 
\item 
$G_{(18,9r)}$ where $r=1,\omega,\omega^2, \omega^3$ when $p\equiv 1\bmod 4$, and $r=1, \nu$ when $p\equiv 3\bmod 4$; 
\item 
$G_{(21,7rs)}$ for $r=0,1,\ldots,p-1$ and $s=1,2,\ldots,(p-1)/2$. 
\end{itemize}
The lists of presentations are generated by taking these parameters 
in a fixed order.

The family $\Phi_1$ consists of the 11 abelian groups;
these are not listed in Table \ref{present}.
For the families $\Phi_{2},\ldots,\Phi_{10}$ 
the rows are partitioned according to the structure of the centre.
Each presentation has a generating set consisting of $\alpha$s and $\beta$s
where the $\beta$s are a minimal generating set for the centre of the group. 
The relation set consists of power relations, commutator relations 
omitting trivial commutators, and relations
expressing central $\alpha s$ in terms of $\beta s$. 
Again there may be one or two additional parameters.

Following \cite{Hall32},
Easterfield \cite{Easterfield40} and James \cite{James80} recorded 
invariants which describe the power structure of 
a regular \mbox{$p$-group} $G$.
Let $\mho^i (G) = \langle g^{p^i} \; : g \in G \rangle$ for $i \geq 0$.
If $|\mho^{i- 1}(G) / \mho^{i}(G) | = p^{w_i}$, 
then define $m_j$ to be the number of $w_i$ which are at least $j$. 
The {\it order type} of $G$ is the partition 
$m_1 m_2 \ldots m_t$ of $6$; 
if $m_{j} = \ldots = m_{j+\delta - 1} = m$ and $\delta > 1$,
then we write $m^{\delta}$. 

The groups of order $p^6$ for $p \in \{2,3,5\}$ are available 
as part of {\sc SmallGroups}. 
While the presentations recorded here 
apply to $p \geq 7$, they also hold for $p = 5$ except for
some groups in $\Phi_{i}$ for $i \in \{35,\ldots,39\}$.

\renewcommand{\baselinestretch}{1.45}

\begin{scriptsize}	
\begin{longtable}[c]{|l|l|l|} 
\caption{Presentations for the groups of order $p^6$ for $p\geq 7$} 
\label{present}\\

\hline
Group & Order type & Presentation \\
\hline
\endfirsthead

\hline
\multicolumn{3}{|c|}{Continuation of Table \ref{present}}\\
\hline
Group & Order type & Presentation \\
\hline
\endhead

\hline
\endfoot

\hline
\endlastfoot


\multicolumn{3}{|c|}{\footnotesize {\bf Groups in $\Phi_{2}$} (Centre of type $C_{p^4}$)}\\
\hline
$G_{(2,1)}$ & $41^{2}$ &  \vtop{\hbox{\strut $ \langle \alpha_{1}, \alpha_{2}, \alpha_{3}, \beta_{1} ~|~ [\alpha_{2}, \alpha_{3}] = \alpha_{1} = \beta_{1}^{p^{3}}, \beta_{1}^{p^4}=\alpha_{2}^{p} =\alpha_{3}^{p} =1 \rangle $} } \\
\hline
$G_{(2,2)}$ & $51$ &  \vtop{\hbox{\strut $ \langle \alpha_{1}, \alpha_{2}, \alpha_{3}, \beta_{1} ~|~ [\alpha_{2}, \alpha_{3}] = \alpha_{1} = \beta_{1}^{p^{3}}, \alpha_{2}^{p} = \beta_{1}, \beta_{1}^{p^4} =\alpha_{3}^{p} =1 \rangle $} } \\
\hline
\multicolumn{3}{|c|}{\footnotesize {\bf Groups in $\Phi_{2}$} 
(Centre of type $C_{p^3}\times C_{p}$ and $\alpha_{1}\in \mho^{2}(Z(G))$)}\\
\hline
$G_{(2,3)}$ & $31^3$ &  \vtop{\hbox{\strut $ \langle \alpha_{1}, \alpha_{2}, \alpha_{3}, \beta_{1}, \beta_{2} ~|~ [\alpha_{2}, \alpha_{3}] = \alpha_{1} = \beta_{1}^{p^{2}}, \beta_{1}^{p^3} = \alpha_{2}^{p} = \alpha_{3}^{p} = \beta_{2}^{p} = 1 \rangle $} } \\
\hline
$G_{(2,4)}$ & $321$ &  \vtop{\hbox{\strut $ \langle \alpha_{1}, \alpha_{2}, \alpha_{3}, \beta_{1}, \beta_{2} ~|~ [\alpha_{2}, \alpha_{3}] = \alpha_{1} = \beta_{1}^{p^{2}}, \alpha_{2}^{p} = \beta_{2}, \beta_{1}^{p^3} =  \alpha_{3}^{p} = \beta_{2}^{p} = 1 \rangle $} } \\
\hline
$G_{(2,5)}$ & $41^2$ &  \vtop{\hbox{\strut $ \langle \alpha_{1}, \alpha_{2}, \alpha_{3}, \beta_{1}, \beta_{2} ~|~ [\alpha_{2}, \alpha_{3}] = \alpha_{1} = \beta_{1}^{p^{2}}, \alpha_{2}^{p} = \beta_{1}, \beta_{1}^{p^3} =  \alpha_{3}^{p} = \beta_{2}^{p} = 1 \rangle $} } \\
\hline
\end{longtable}
\end{scriptsize}

\bigskip
\noindent
\begin{tabbing}
456789012345123454545555555545\=4555555555555555555555555545\=
4555555555555555555555555545 \kill
Mathematical Sciences Institute \>  Department of Mathematics \> 
Christ Church \\
Australian National University \> University of Auckland  \>  
University of Oxford \\
ACT 0200 \> Auckland  \> OX1 1DP\\
Australia \> New Zealand \> United Kingdom \\ 
mike.newman@anu.edu.au \> e.obrien@auckland.ac.nz  \> 
michael.vaughan-lee@christ-church.oxford.ac.uk
\end{tabbing}

\begin{thebibliography}{99}

\bibitem{Small} Hans Ulrich Besche, Bettina Eick and E.A. O'Brien. 
\newblock A millennium project: constructing small groups.
{\it Internat.\ J.\ Algebra Comput.}, {\bf 12}, 623--644, 2002.

\bibitem{Magma}
Wieb Bosma, John Cannon, and Catherine Playoust.
\newblock The {\sc Magma} algebra system I: The user language.
\newblock {\em J.\ Symbolic Comput.}, {\bf 24}, 235--265, 1997.

\bibitem{Easterfield40}
Thomas~E.\ Easterfield.
\newblock A classification of groups of order $p^6$.
\newblock PhD thesis, Cambridge University, 1940.

\bibitem{GAP4} The GAP Group, GAP -- Groups, Algorithms, and Programming, Version 4.12.1; 2022, \href{https://www.gap-system.org}{www.gap-system.org}
See also \href{https://docs.gap-system.org/pkg/smallgrp/doc/chap1.html}
{docs.gap-system.org/pkg/smallgrp/doc/chap1.html}

\bibitem{Hall32}
P.\ Hall. 
\newblock {A contribution to the theory of groups of prime-power order.}
\newblock {\em Proc.\ London Math.\ Soc.}, {\bf s2-36}, 29--95, 1934.


\bibitem{James80}
Rodney James.
\newblock The groups of order $p^6$ ($p$ an odd prime).
\newblock {\em Math.\ Comp.}, {\bf 34}, 613--637, 1980.


\bibitem{p6}
M.F. Newman, E.A. O'Brien, and M.R. Vaughan-Lee,
\newblock Groups and nilpotent Lie rings whose order
is the sixth power of a prime,
{\it J. Algebra}, {\bf 278}, 383-401, 2004.

\bibitem{github}
M.F. Newman, E.A. O'Brien and M.R. Vaughan-Lee.
The groups of order $p^6$. 
\href{https://github.com/eamonnaobrien/The-groups-of-order-p-6/}
{github.com/eamonnaobrien/The-groups-of-order-p-6}

\bibitem{error-list}
E.A.\ O'Brien, Sunil Kumar Prajapati, Ayush Udeep.
Comment on: The groups of order $p^6$ ($p$ an odd prime) 
By Rodney James, {\it Math.\ Comput.\ }{\bf 34} (1980), 613--637.
\href{http://arxiv.org/abs/2308.04443}{arxiv.org/abs/2308.04443}

\bibitem{mrvl}
Michael Vaughan-Lee.
The nilpotent Lie rings of order $p^k$ for $k \leq 7$.
\href{http://www.iaa.tu-bs.de/beick/soft/liepring/p567.pdf}
{www.iaa.tu-bs.de/beick/soft/liepring/p567.pdf}

\end{thebibliography}
\end{document}